\documentclass[pdflatex,sn-mathphys-num]{sn-jnl}%

\usepackage{graphicx}%
\usepackage{multirow}%
\usepackage{amsmath,amssymb,amsfonts}%
\usepackage{amsthm}%
\usepackage{mathrsfs}%
\usepackage[title]{appendix}%
\usepackage{xcolor}%
\usepackage{textcomp}%
\usepackage{manyfoot}%
\usepackage{booktabs}%
\usepackage[font=small,skip=0pt]{caption}
\usepackage[font=small,skip=0pt]{subcaption}
\usepackage{calligra}
\usepackage{pbsi}
\usepackage{tikz}
\usepackage{array}
\usepackage{dashrule}
\usetikzlibrary{arrows.meta, positioning}
\usetikzlibrary{math}
\usetikzlibrary{lindenmayersystems}
\usetikzlibrary{decorations.fractals}

\newcommand{\letterbreak}{%
    \par\bigskip
    \noindent\rule{\textwidth}{0.4pt}
    \bigskip
}

\newcolumntype{L}[1]{>{\raggedright\arraybackslash}p{#1}}

\begin{document}

\title[]{Letter to a Young Black Mathematician: History, Storytelling, and Belonging in Mathematics}


\author[1]{%
  \fnm{Ranthony A.} \sur{Clark}
}

\affil[1]{\orgdiv{Department of Mathematics}, \orgname{The Ohio State University}}

\abstract{What does it mean to tell the story of a mathematician? Written as a letter to a young Black mathematician, this article reflects on how stories shape our understanding of mathematics and our sense of belonging within the mathematical community. Through the lens of Hidden Figures Revealed, a project documenting Black mathematicians at The Ohio State University, I discuss humanity, mathematics, and context as three guiding principles for mathematical storytelling. In particular, I reflect on how this framework emerged through archival research, shaped historical narratives, and informed mathematical learning experiences.} 


\maketitle


\section*{Dear Young Black Mathematician,}

\begin{quote}
\emph{``You're the first Black woman I've seen work in the department in the last twenty years."} 
\end{quote}

I was in the break room. At this point I don't remember why anymore. I could have been warming up lunch or picking up a stack of papers from the printer; I just remember her voice. It was gentle, but husky. She sat among the other Black custodial workers. They were trying to think if they could recall another Black woman during their years cleaning the building on Ohio State's campus that housed the math faculty and postdocs. They couldn't. I remember them nodding in approval as if my being there was something noteworthy to witness. I remember the warmth in her words and in their presence, and how the energy of that moment was captured over the several years in smiles across rooms, nods in elevators, and casual chatter during the times I crossed paths with them throughout the department. I remember feeling that they saw me, and during those moments when I saw them, I was not alone.

My first postdoc brought feelings of isolation I hadn't experienced coming from a graduate program where I had never been the `only' across a variety of axes of diversity. And yet I found myself entering a department where I was the only Black woman among nearly two hundred faculty, staff, postdocs, and graduate students. I wonder if you know this feeling. A haunting awareness that your aptitude for this subject facilitates your belonging, but that parts of who you are could also keep you apart from its community. The strange experience of loving mathematics, yet sensing you are standing on the outside while being inside: a kind of mathematical double consciousness:\footnote{\emph{The Souls of Black Folk} \cite{dubois1903souls} is a seminal text by W.~E.~B. Du Bois that develops the notion of \emph{double consciousness} which he notes as the experience of Black Americans seeing themselves simultaneously through one's own eyes and through the eyes of a society that understands Blackness as other.'' I use the phrase ``mathematical double consciousness'' here to describe the experience of seeing oneself simultaneously as a Black person and as a mathematician. Note that I speak of my own experiences here, and do not presume that every Black mathematician feels tension between these dual perspectives.} ``\emph{one ever feels his twoness}"  \cite{dubois1903souls}.

I wish I could say that I took those feelings of isolation and paid them no mind as I navigated the research and teaching milestones that came at the beginning of my postdoc. I am certain that from the outside I looked well adjusted: I was doing ``all the things.'' But there were times I was stifled and could not get out of my head. The seminar I sat in silently for an entire semester, too afraid to speak. What if I was wrong? What if my question made no sense? The moments I couldn't immediately answer a student's question in lecture, the insecurity welling up deep in my stomach was enough to make me dry heave after class. I was representing something bigger than myself here. I was the only Black woman. I could not mess up.

There were other moments. The dinner with the speaker where I talked about my Hidden Figures course, sharing that we taught the women at NASA were working to get to the moon, while in Virginia, in the same era, there was a generation of Black children denied education, prompting a tangent from my postdoc mentor at the time about misunderstood historical governors.\footnote{\emph{Context}: At our seminar dinner we discussed teaching assignments and I was describing the Hidden Figures course I co-designed. I shared it was remarkable how much progress was made within NASA by the Black and white women computers, in contrast with the segregated conditions within Virginia, prompting a brief discussion about how the actions of some historical governors were misunderstood. In September 1958, Governor J. Lindsay Almond Jr. ordered public schools in Charlottesville, Front Royal, and Norfolk closed rather than allow desegregation. In Prince Edward County, officials closed public schools in 1959 rather than comply with a federal desegregation order. Although the county's public schools were closed, many white families were able to establish and attend segregated private schools, while most Black children were left without meaningful access to schooling. The public schools remained closed until 1964 \cite{bonastia2009white}.} While he talked about how George Wallace was frequently misjudged, I texted under the table to my father who grew up in Birmingham, Alabama in the 1960s, `\emph{was George Wallace the one that said segregation segregation segregation?}' The phone buzzed with an incoming message with a slight correction, `\emph{segregation now, segregation tomorrow, segregation forever}\cite{wallace1963inaugural}.' I tried to write off this discussion as an intellectual exercise, but then there was the research meeting with the same mentor that started casually with the remark, `\emph{you know, when people see you, they see diversity.}' I responded I could not help how people saw me, which seemed to satisfy him, but I was too frazzled to focus afterwards. 

Still, I kept moving. Producing. Stuffing down the feelings of isolation and otherness as neatly and meticulously as I could. In the summer of 2020, when the \#ShutDownSTEM initiative called on STEM researchers to pause their usual activities for a day of reflection and action around antiBlack racism in academia and STEM, those feelings came rushing back. I hope you haven't felt the need to charge on like this, to avoid the discomfort. I wish I had admitted to myself sooner it was not sustainable.

I was two years into my postdoc then. I found myself reflecting on what Phyllis had said. I remembered her voice and the sentiment, \emph{you're the first Black woman here in 20 years}. Was it possible that I was the only Black woman to ever work for the department? I had to know. While the world was `shutdown,' I began poring over the history of the department where I worked, desperate to see if what I believed was true, that there were other Black mathematicians who preceded me, that I was not alone.

\begin{figure}[t]
    \centering
    \includegraphics[width=\textwidth]{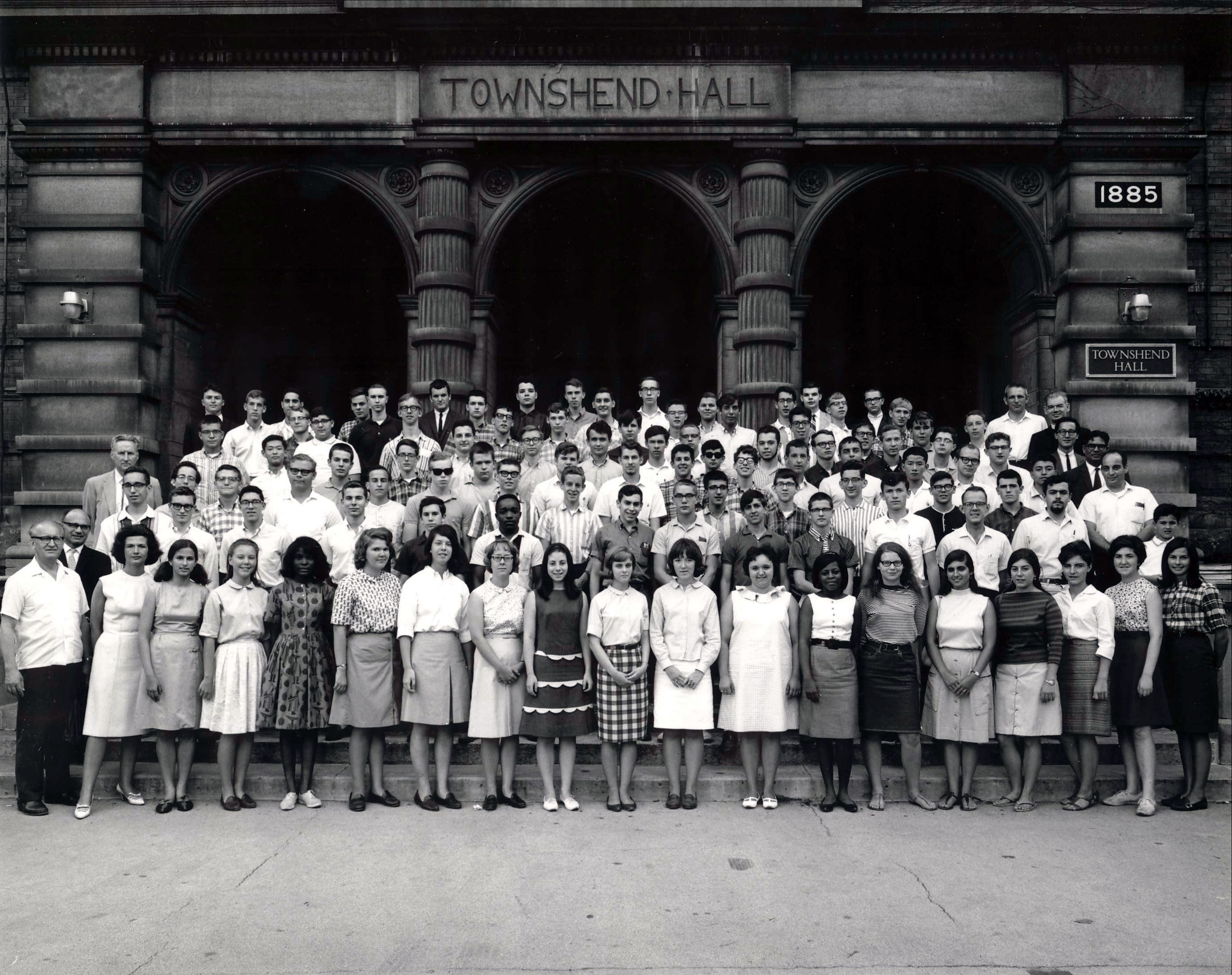}
    \caption{Rada Higgins McCreadie (featured in first row, fifth from left in a polka dot dress) as a participant in the Ross Mathematics Program at The Ohio State University, summer 1965.}
    \label{fig:rada_ross}
\end{figure}

This inclination is what led me to her. In 1974 Rada Higgins earned a PhD in mathematics from Ohio State. She was later identified as the 20th Black woman in the United States with this degree. She was my existence proof. I found an old black and white photo of her in the archives for the 1965 Ross Mathematics Program at The Ohio State University. There she was, studying elementary number theory under Arnold Ross, sharp in polka dots in the front row of her peers, arms tucked slightly behind her back (Figure~\ref{fig:rada_ross}). I wonder if she spoke up in groups. I wonder if she doubted herself when she was asked questions. I wonder if she felt like she belonged. She applied to Ohio State for graduate school and began her doctoral studies in 1969. Over 50 years later I was starting in the same department. I had a named postdoc as a Ross Assistant Professor, in honor of the mentor who supported her.

I had gone looking for proof that someone like me had existed in the department before me, and I found her. It was like I was an undergraduate again at my first math conference, running into the first Black woman mathematician I had ever seen, Jana Talley. I asked her for a hug, asked if she really had a PhD. Something in me changed after meeting her, her very existence shifted a PhD from something that felt overwhelming to something I could possibly do. Like she noted, I just needed to find support. I did not interact with Rada. We did not embrace. But knowing that I was not the first changed me. It lifted the burden I felt from being the ``only.'' I wanted to know more. I wanted to know about her research, her family, and her career. I wanted to know what it was like for her, what it meant to become a Black mathematician in Columbus in the 1960s and 1970s. I wanted to make her existence known, and to learn more about other Black mathematicians from the department.

This became a personal fixation. It felt important that whoever came after me would not be faced with some version of \emph{you're one of a handful of Black people to matriculate through here in the last two decades.} I reached out to departmental personnel for historical records and scoured the internet for evidence, mining sites like the Mathematicians of the African Diaspora (MAD) \cite{williamsMAD} and Mathematicians of the African Diaspora Database's Ensemble of Researchers (MADDER) \cite{madder, walker2019existence} to find as many Black people from Ohio State as I could who had earned an undergraduate, master's, or doctoral degree in mathematics.\footnote{The \emph{Mathematicians of the African Diaspora} (MAD) is an online repository documenting Black mathematicians and their contributions \cite{williamsMAD}. It was created and maintained by Scott Williams until his retirement. The \emph{Mathematicians of the African Diaspora Database's
Ensemble of Researchers} (MADDER) is ongoing and expanded Williams' work into a modern searchable database of Black mathematicians \cite{madder, walker2019existence}.} In time, I learned that nearly 200 mathematicians who identified as Black had earned degrees from The Ohio State University.\footnote{This number was confirmed through the Hidden Figures Revealed project described in \cite{clark2025hiddenfigures}.} They went on to careers spanning academia, industry, and education: university presidents, authors, lawyers, and mathematicians. And yet to many, their stories remained hidden.

That search ultimately became the project ``Hidden Figures Revealed: Dynamic Narratives and History of Black Mathematicians from The Ohio State University." With the support of collaborators and grant funding, what began as a solitary archival effort grew into the first comprehensive historical study of Black mathematicians at a single institution \cite{clark2025hiddenfigures}. But here I want to be honest with you: this project started with a need to affirm, first for myself, then for you, that we are not alone. As I learned more, I began to ask other questions: What does it mean to tell the story of a mathematician? What gets preserved and what disappears? What happens when we refuse to choose between the person and the mathematics? What if a mathematical story could attend to the person, the mathematics, and the world in which that mathematics took place?

There are many ways to respond to stories and history. There are many ways to work toward a mathematical community in which more people can see themselves belonging. I am writing to you now from the beginning of my journey as a research professor, wondering what I can do to make my own way while clearing a path for you. One day, you too may begin to think about if, when, or how you want to do something concrete to make mathematics more inclusive for yourself or for others. As for me, I choose history. I choose storytelling.\\


\letterbreak


\begin{quote}
\emph{``I'm interested in the way in which the past affects the present and I think that if we understand a good deal more about history, we automatically understand a great deal about contemporary life.''}

\hfill -- Toni Morrison \cite{morrison2008what} --
\end{quote}

As you learn mathematics, you will encounter many names. Some will come from theorems, some from techniques, some from conjectures, and some from entire areas of mathematics. If your experience has been anything like mine, there may not always be a story attached to the mathematics you encounter. I learned about math without mathematicians. Even some of my favorite subjects in school were presented in ways that were stuffy, formulaic, or abstract, with very little attention to the people who developed the mathematics or to the political and social contexts in which they lived. The more historical work I have done, the more I have come to see how much the past can tell us about the present. This is certainly true in mathematics. Knowing where an idea came from, who developed it, what problem they were trying to solve, and what was happening around them has often made the mathematics feel more alive to me. 

I want to encourage you to look for mathematical stories. I will share one of my favorites. When I was younger, I learned the quadratic formula through a song set to the tune of ``Pop Goes the Weasel." It served as a mnemonic to memorize: 
$$\displaystyle x = \dfrac{-b \pm \sqrt{b^2-4ac}}{2a} $$
\noindent where ``ALL over 2 - A'' was analogous with the ``POP goes the weasel'' in the original song. Years later as an undergraduate, I read William Dunham's \emph{Journey to Genius}, and learned there were formulas for roots of cubic and quartic polynomials! I learned about the secrecy, the oaths, and the rivalries surrounding the discovery and publication of solutions to the cubic among Cardano, del Ferro, Tartaglia, and Lodovico Ferrari \cite{dunham1990journey}. I also learned that there was no such formula for quintic polynomials.  

I was wildly entertained by the tragic tale of Évariste Galois and his work on whether there was a formula for solving quintic polynomials. The answer was no, and he was able to prove it using the symmetries of the roots of polynomials~\cite{stewart2003galois}. He sadly died in a duel, infamously spending the night before writing up notes on his mathematical ideas~\cite{rothman1982short}. It took 14 years for his manuscripts to be published, and another century before the theory was formalized in the way that we learn it today in abstract algebra~\cite{Galois1846, artin1998galois}. The immortalization of his work often tidies up the struggles he faced in life. He tried to send his work to Cauchy and Poisson. Rejected~\cite{rothman1982short}. He enrolled in what is now École Normale Supérieure to study mathematics. Expelled~\cite{toti1998evariste}. He frequently pushed back against the French monarchy. Jailed~\cite{toti1998evariste}. 
 
I don't know where all of this went when I was introduced to group theory. By then, Galois was a name attached to theory. His history had been separated from the mathematics. I understand why this happens. Especially now as a professor, I understand the constraints and pressure to progress through the curricula. But the ultimate consequence of my training is that, aside from isolated instances, I learned to think about mathematics independently from mathematicians. Is it any wonder that I was 22 years old when I realized other Black women mathematicians existed, or that I had no idea Ohio State had a rich history in cultivating Black mathematical talent?  

In graduate school I got the book \emph{Hidden Figures: The American Dream and the Untold Story of the Black Women Who Helped Win the Space Race} by Margot Shetterly~\cite{shetterly2016hidden} as a gift after passing my comprehensive exam to advance to candidacy. I could not get enough of the story. I cried when I saw the movie. I imagine you are familiar with the story. I relished the attention it received, but there were times I felt uneasy about the framing of these women as `\emph{sentimental firsts}' \footnote{In \emph{Between the World in Me}, Coates laments history books that ``spoke of black people only as sentimental ‘firsts’—first black five-star general, first black congressman, first black mayor—always presented in the bemused manner of a category of Trivial Pursuit." Here he encourages a treatment of Black history that is just as serious as other Western histories and not reduced to heroism, struggle, or novelty \cite{coates2015between}. I was worried that many were treating the women from the Hidden Figures story this way.} \cite{coates2015between}. I was excited they were being celebrated, but I wanted to know, what mathematics did they actually do? What problems were they solving? How were their opportunities shaped by the historical and political circumstances surrounding the Space Race? 

I developed a Hidden Figures course based on their story with John Johnson to try and answer these questions by exploring the mathematics of the Space Race through the eyes of its human computers~\cite{osu2019hiddenfigures, edmonds2021intersections}. We purchased slide rules from the International Slide Rule Museum to lend authenticity to the daily tools they might have used in their work. (I developed an affinity for the Post 1460 Versalog.) We also examined the historical context influencing their opportunities in mathematics in the mid-20th century, as well as their communities, personal lives, and backgrounds. 
As I worked on the course, I realized this was no different from how I felt about the stories of Cardano, Tartaglia, Ferrari, and Galois. I didn't want to choose between the person and the mathematics. The story made me crave the math.

When I worked on \emph{Hidden Figures Revealed}, I knew I wanted to take a similar approach. I did not want to reduce people we learned about to their mathematical contributions, and I did not want their mathematics to disappear beneath their biographies. I wanted to know about the person, their mathematics, and the world in which they were doing that mathematics. That became the framework I kept returning to. 

\begin{figure}[htpb]
\begin{center}
\begin{tikzpicture}[
    every node/.style={font=\small},
    mainnode/.style={
        circle,
        draw=black,
        thick,
        minimum size=2.5cm,
        align=center,
        fill=gray!10
    },
    arrow/.style={
        <->,
        thick,
        >=Stealth,
        rounded corners=8pt
    }
]

\def\radius{2.25}

\node[mainnode] (A) at (90:\radius) {Humanity};
\node[mainnode] (B) at (210:\radius) {Mathematics};
\node[mainnode] (C) at (330:\radius) {Context};

\draw[arrow] (A) to[bend right=25] (B);
\draw[arrow] (B) to[bend right=25] (C);
\draw[arrow] (C) to[bend right=25] (A);

\end{tikzpicture}
\end{center}
\caption{Three principles for mathematical storytelling: humanity, mathematics, and context.}
\label{fig:storytelling_principles}
\end{figure}
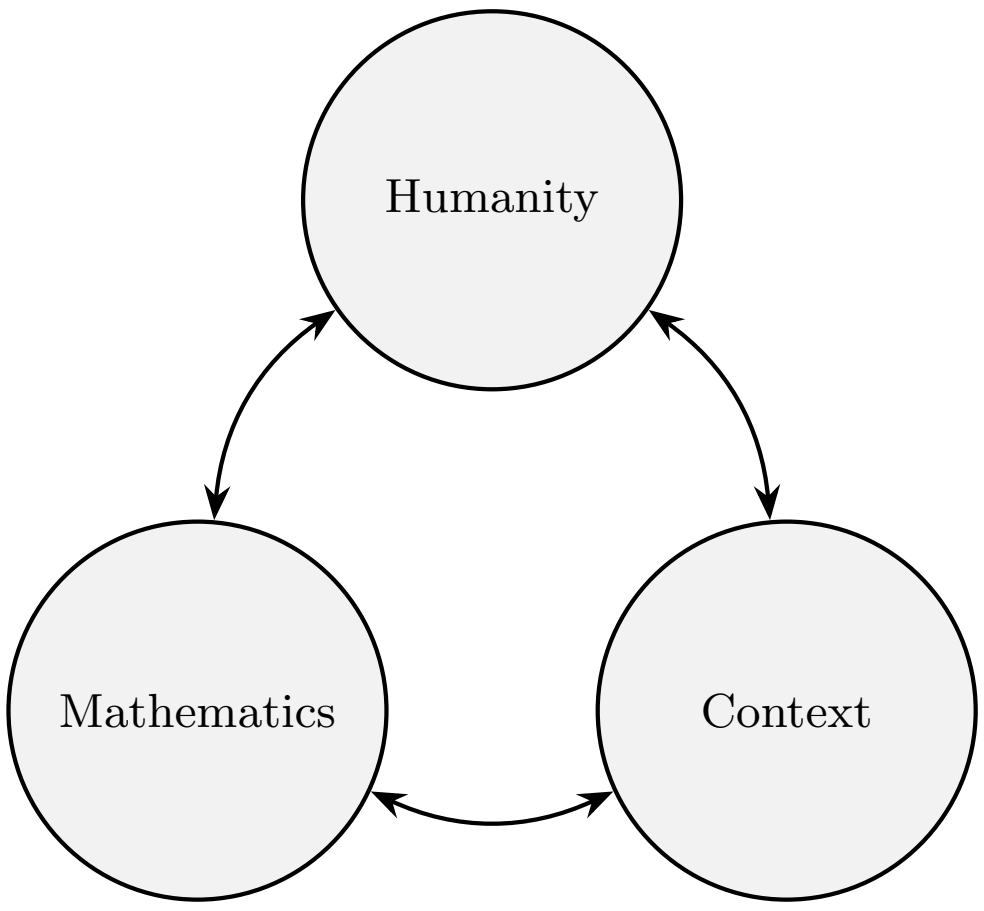

I came to think of these as three principles for mathematical storytelling: \textbf{humanity, mathematics, and context} (Figure~\ref{fig:storytelling_principles}). I don't think there is a hard and fast rule about the weight each should carry in a particular story. Some stories may lean heavily toward the mathematics. Others may require more historical or biographical context. But I don't think any of the three should disappear entirely.

These principles guard against three different kinds of incomplete stories. Humanity without mathematics can turn a mathematician into an inspirational symbol whose actual intellectual work disappears. Mathematics without humanity can make that work feel disembodied, as though ideas emerge independently of the people who develop them. And both can become misleading without context. A career can begin to look inevitable, or like the product of individual talent alone, when it was also shaped by institutions, communities, historical events, and opportunities. Context without mathematical content creates a different problem. It may tell us something important about a person's life without telling us much about their life as a mathematician. For me, the point is not to give these three dimensions equal weight, but to resist allowing any one of them to erase the others.

 There are, of course, longstanding traditions for how to do this that don't explicitly invoke these three principles. I am not trying to invent a new method of storytelling. Scientific biographies have long connected people's lives to the development of scientific ideas~\cite{soderqvist2007history}, and historians of science have reminded us that scientific work does not happen in isolation from the intellectual and social worlds around it~\cite{shapin1994social}. Literary journalism offers another tradition of paying close attention to people, events, and context~\cite{kramer1995literary}. You may remember learning that a good story should answer the questions `who, what, where, when, why, and how.' These questions help us establish what happened. But when we tell the story of a mathematician's life, there are other questions worth asking, too.
 What do we choose to include? What do we leave out? What do we make visible about a person's life, and what disappears from view?

People who study mathematics education have asked some of these questions by using mathematical stories and autobiographies to understand how people experience mathematics and how those experiences become part of how they see themselves as mathematicians, teachers, or learners~\cite{drake2006turning,mcculloch2013math}. This matters to me because I believe the stories we encounter can shape what we think a mathematician looks like, what a mathematical life can look like, and whether we can imagine ourselves in one.

So how do you tell a multifaceted mathematical story? When you tell stories about mathematicians, incorporate humanity by asking the question: \emph{Who was this person?} There are the classics. Where someone went to school, what degree they earned, or where they worked. But you should dig beyond pedigree whenever possible.  Where do they come from, what did they care about, what did they do outside of mathematics, what communities were they a part of, and what experiences shaped their pathway into and within mathematics? 

The stories we tell should also answer the question: \emph{What mathematics did they do?} This one seems obvious, but I think it is worth saying. If we are telling stories about mathematicians, I don't want the mathematics to become incidental. There should mathematics in a story about a mathematician. This presents many challenges due to vastness of our subject and the myriad of specialties. But it also presents a neat opportunity to dive into the ideas that lit somebody else up, to immerse yourself in familiar or completely new research areas. I want people to encounter the ideas that these mathematicians worked on, even if that means finding creative ways to make those ideas accessible for themselves or others. That could mean narrowing to one aspect of a person's research program, like a dissertation or a single paper. It could mean developing an activity that allows a student to encounter their mathematical ideas without having to understand the mathematics at the level at which it was originally developed. Pick one place to start. One idea to digest. From there you can decide the level of breadth and depth of the mathematical story you want to share.  

We should also consider the question: \emph{What was happening around them?} Mathematical lives don't happen in a vacuum. They navigate shifts in colloquial knowledge as well as competing frameworks for moving central ideas in their fields forward. They live through political movements and social changes. They encounter barriers, opportunities, institutions, and communities. Mathematical stories should include not just who someone was or the mathematics they did, but what was happening around them that helped shape the path they took. How can we paint the political, social, and historical backdrop that colored their experiences within the mathematical community? 

There is a responsibility that comes with all of this. Historical records are incomplete. We may sometimes know what someone did but not why they did it, or what happened in time but not how someone experienced it. We can't project our imaginings. We can't make things tidy because it is convenient. Not every person has to become a hero, a symbol, or a cautionary tale. We have to lead with honesty and be sincere about what we know, what we infer, and what we can't corroborate. 

Even if you are not called to tell the stories of other mathematicians, I encourage you to think about these principles to process your own mathematical journey. The part that usually takes the most time when I process my own story is context. I was not always used to thinking about my research as influenced by external events.

Here's an example. In 2020, there were widespread racial justice protests across the world after the murder of George Floyd. The \#ShutDownSTEM initiative was a reaction to this. Locally, Ohio State established a seed fund for racial justice, and I earned one of the inaugural awards to lead a team to study the Black math history of the department. Meanwhile, insecure about my relationship with my research mentor and more broadly disillusioned with my doctoral research,  I started to work with a group in the department under a Research Training Group\footnote{The National Science Foundation's funds Research Training Groups that foster collaborative research environments that enhance training for graduate students and postdocs. The group at Ohio State was supported by NSF \#DMS-1547357 RTG: Algebraic Topology and Its Applications.} that studied Applied Algebraic Topology. At the same time, during the 2021 redistricting, Moon Duchin's lab was collaborating with organizations in Columbus as well as Ohio State to do community mapping to push for fairer maps from the state's new legislative redistricting commission~\cite{chambers2021ohio,mggg2021cois}. She reached out to one of my mentors in the department for help. He thought it would be a good opportunity and recommended me.

Then there was another piece of context I did not fully appreciate at the time. On his first day in office, President Joe Biden signed an executive order to advance support for racial minorities in the federal government~\cite{biden2021executive}. Federal agencies developed strategic plans in response~\cite{nsf2021equity}. Part of the NSF's plan included a new postdoctoral fellowship program designed to broaden participation in the sciences~\cite{nsf2021mpsascend}. I wrote a proposal for this new award on how to think about applications of metric geometry and applied algebraic topology to electoral redistricting. I emphasized my Hidden Figures work in my broader impacts. I was an inaugural recipient of that award. It changed my life and my research program. 

There is a version of this story where I simply say I changed research directions and received a fellowship. It's technically true, but it's too tidy. The actual story is full of people, institutions, politics, funding opportunities, and some serendipity. My mathematical trajectory did not happen in isolation from the world around me. Neither did the trajectories of the mathematicians I had been studying. This framework helped me take a list of names of Black mathematicians from Ohio State and begin to see the people, mathematics, and worlds behind those names. I came to see that the three principles aren't independent. The person influences the mathematics, the mathematics interacts with the world around them, and the world around them shapes what mathematical opportunities are available to a person. This way of learning and telling stories has helped me feel deeply connected to mathematics, even when I don't always see myself in the broader mathematical community. If this is something you have ever struggled with, maybe it can help you too.\\


\letterbreak


\begin{quote}
\emph{In any category of ``firsts," there is always an underlying question of how ``the first" decided that it was within his or her power to become that which no one like him or her had become before. The answers Elbert Frank Cox and Euphemia Lofton Haynes (the first Black male and female Americans to earn doctoral degrees in mathematics) might have given to this question, unfortunately, are unknown. But what we do know about them--largely through remembrances from Cox's students, colleagues, and family members and from Lofton Hayne's collection of family and professional papers held at Catholic University--tells us something about their journeys as mathematicians and, also, the ways in which these journeys left blueprints for others to follow.}\\

\hfill -- Erica Walker, \emph{Beyond Banneker} \cite{walker2014beyond} --
\end{quote}

I remember lighting a candle and taking a deep breath to stay calm as I logged into the Zoom link she had sent me. She was squeezing in time during her office hours, because I wanted to get some advice on how to approach my interview protocol. The project and the grant that supported it had grown. What began as my attempt to find the Black mathematicians who came before me had become a research project with a grant, a team, and plans to hire student fellows. I was now responsible for leading people through a project that I was still learning how to do myself. 

I wish I could tell you that you will always feel ready for the things you are asked to do. You won't. Sometimes you will be uncomfortable because you are doing something genuinely new, and I hope you learn not to mistake that discomfort for evidence that cannot do it and that you don't belong. Ask for help. Find someone who knows what you don't know. I was the least senior person on our grant's leadership team. In fact, I had to get special approval from my department to be a principal investigator as a postdoc. Our study had to be approved by the Institutional Review Board.\footnote{Institutional Review Boards review research involving human participants, including interviews, to consider questions of informed consent, privacy, confidentiality, and potential risk to participants. For me, this was one of the first indications that historical research of living mathematicians involved forms of responsibility that were quite different from those I had encountered in mathematical research.} This was completely different than mathematical research, and even that research was new to me as I moved into a new field. Everything was new. I was overwhelmed. I felt the familiar insecurity welling up deep in my stomach. I could not mess up. I wanted to lead well. I wanted to tell their stories right. But how? 

Erica Walker was one of few people I knew who had done something similar. I read \emph{Beyond Banneker: Black Mathematicians and the Paths to Excellence} \cite{walker2014beyond} in 2020. I can't remember if it was before or after \#ShutDownSTEM, but I knew that when we started the Hidden Figures Revealed project I wanted to use it as a guide. Our meeting and her book gave me a model for how interviews, historical context, and individual stories could be woven together without reducing the mathematicians to their identities or their achievements. She had interviewed thirty five mathematicians and organized their stories into three time periods, each marked by common historical experiences, including legal and de facto segregation. We adopted a similar approach for Ohio State, using the degrees earned by Black mathematicians in the department to help identify three research periods. The seven Black mathematicians who earned PhDs in mathematics from Ohio State all received their degrees between 1963 and 1984. We called them the Pioneers and this became Research Period 1. This was the period we studied the most.

Have you ever considered how intimately the Civil Rights Movement and other fights for racial equity in the United States are connected to mathematics? I had not before this project. In 1979, five years after Rada Higgins finished her PhD, Columbus was still dealing with the consequences of the 1954 \emph{Brown v. Board of Education} Supreme Court ruling. Over twenty years later, the 1977 case \emph{Penick v. Columbus Board of Education} from the U.S. District Court ruled the city still had de facto segregation in its public school system and needed to desegregate. The city fought this ruling up to the Supreme Court and lost. In 1979 they began desegregation via busing. There was a dramatic decline in Columbus public school enrollment in the decades that followed. I cannot say from the evidence we found how these changes affected the Ross Program's local recruitment, but they formed part of the environment in which the program operated. ``White flight'' to the suburbs, along with other families and educators abandoning city schools, led to the closure of several schools. East High School, where Rada Higgins went to school, remains open today. Another school, Central High School, was forced to close its doors due to low enrollment. Part of it has been historically preserved as part of the Center of Science and Industry (COSI). This area was part of my jogging route when I was working on this project. 

The juxtaposition is striking. Carolyn Mahoney's first year of her PhD was in 1979. She was the last Black woman to earn a PhD in mathematics from Ohio State, earning her degree in 1984. She and Rada Higgins were both able to pursue graduate studies at Ohio State with support from the National Science Foundation. They pursued PhDs during a period of increased federal investment in science and mathematics education associated with the Space Race and the National Defense Education Act of 1958. This expansion of federal support helped create new opportunities for mathematical research, education, and enrichment, including programs like the Ross Mathematics Program and fellowship opportunities for faculty at Historically Black Colleges and Universities (HBCUs). Guy Hogan and Thyrsa Frazier Svager were able to take advantage of a special NSF program that provided fellowships to faculty at Historically Black Colleges and Universities (HBCUs) that were teaching but did not yet have a doctoral degree \cite{clark2025hiddenfigures}.  

Everything we learned provided more context. We could not provide a vivid picture of the lives of every single Black math PhD from Ohio State during this time period. We could, however, take care to understand key historical, institutional, social, and political factors that formed the environment in which mathematicians at Ohio State studied and worked during each Research Period. I followed Erica's advice and got support. Support from social scientists at Ohio State to refine our interview protocol. Support from our community partner, the National Afro American Museum and Cultural Center in Wilberforce, Ohio, to connect us to faculty at Central State University, also in Wilberforce, the HBCU where Thyrsa Frazier Svager and Guy Hogan worked before completing their PhDs at Ohio State. 

\begin{figure}[p]
    \begin{center}
    \fbox{\includegraphics[
        width=\textwidth,
        height=\textheight,
        keepaspectratio
    ]{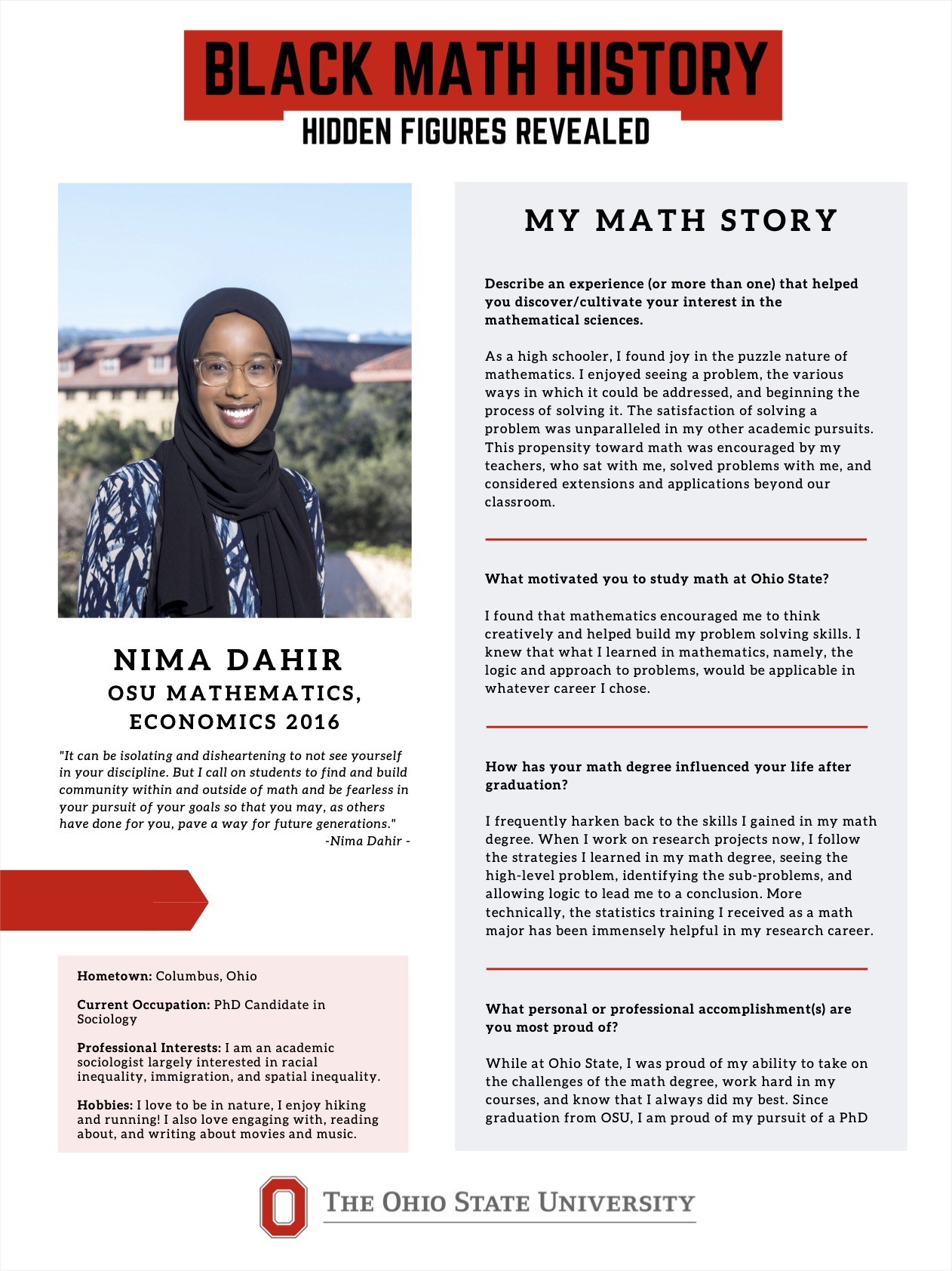}}
    \end{center}
    \caption{Alumni profile of Nima Dahir, now Assistant Professor
    of Sociology at Ohio State.}
    \label{fig:profile_example}
\end{figure}

We took our students to this museum. At the time there was an exhibition \emph{Queens of the Heartland} featuring illustrated portraits with panel text for thirty women who were leaders of social change with deep connections to Ohio. We took notes on the way their stories were told in this simple way, which eventually inspired public facing profiles for Black math alumni\footnote{These profiles were inspired by similar efforts like \emph{Lathisms} and \emph{Mathematically Gifted and Black}.  I built the questions used from the profile format used by \emph{Mathematically Gifted and Black}, with permission. } (Figure~\ref{fig:profile_example}). We hoped that the profiles would simultaneously highlight Black alumni, the humanity of their story, and what sort of careers a math degree could lead to. We gave their archivists names of the Black math PhDs we knew from Ohio State. They directed us to a room where we could look inside, but not take, materials in the Helen McWorter-Simpson Collection.\footnote{The Helen McWorter-Simpson Collection (NAM MSS 2019), acquired in 1993, is held by the National Afro-American Museum and Cultural Center in Wilberforce, Ohio. The finding aid \cite{mcworterSimpsonCollection} describes a 38-box collection spanning ca. 1830--1993, comprising correspondence, family documents and histories, photographs, newspapers, legal notices, U.S. military papers, patents, and materials from various organizations.} It contained newspapers, correspondence, personal papers, patents, and other materials from the McWorter-Simpson families.
Helen was the great granddaughter of `Free Frank' McWorter, the first non Native American settler of Hadley, Pike County, Illinois. Her father was awarded a patent for the invention of the autoplane, a predecessor of the helicopter, her daughter won an Emmy. She was also a close relative of William McWorter, the first Black mathematician to earn a PhD at Ohio State in 1963.

Most archival work constructs simple narratives, when they arrived at Ohio State, what degrees they earned, where they worked, and sometimes what they had accomplished. But it does not always answer the question Erica Walker had asked about the ``firsts'': how does someone come to believe that it is within their power to become something that no one like them had become before? Inside the McWorter-Simpson collection was a physical copy of \emph{Of The Structure of pre-$p$-Rings}, by William McWorter, Jr. I had no way of knowing whether anyone else in the family knew what a ring was, let alone a $p$-ring, or a pre-$p$-ring, but it was there among photographs from weddings and letters. There was something striking about encountering his mathematics not as a citation in a university record, but as a physical object preserved alongside the papers of his family. We could not know from the archive whether McWorter's family history shaped his mathematical imagination. But he came from a family of firsts, inventors, and pioneers. What does it mean to grow up surrounded by evidence that the boundaries of what is possible can be crossed?

\begin{figure}[htpb]
    \begin{center}
    \includegraphics[width=\textwidth]{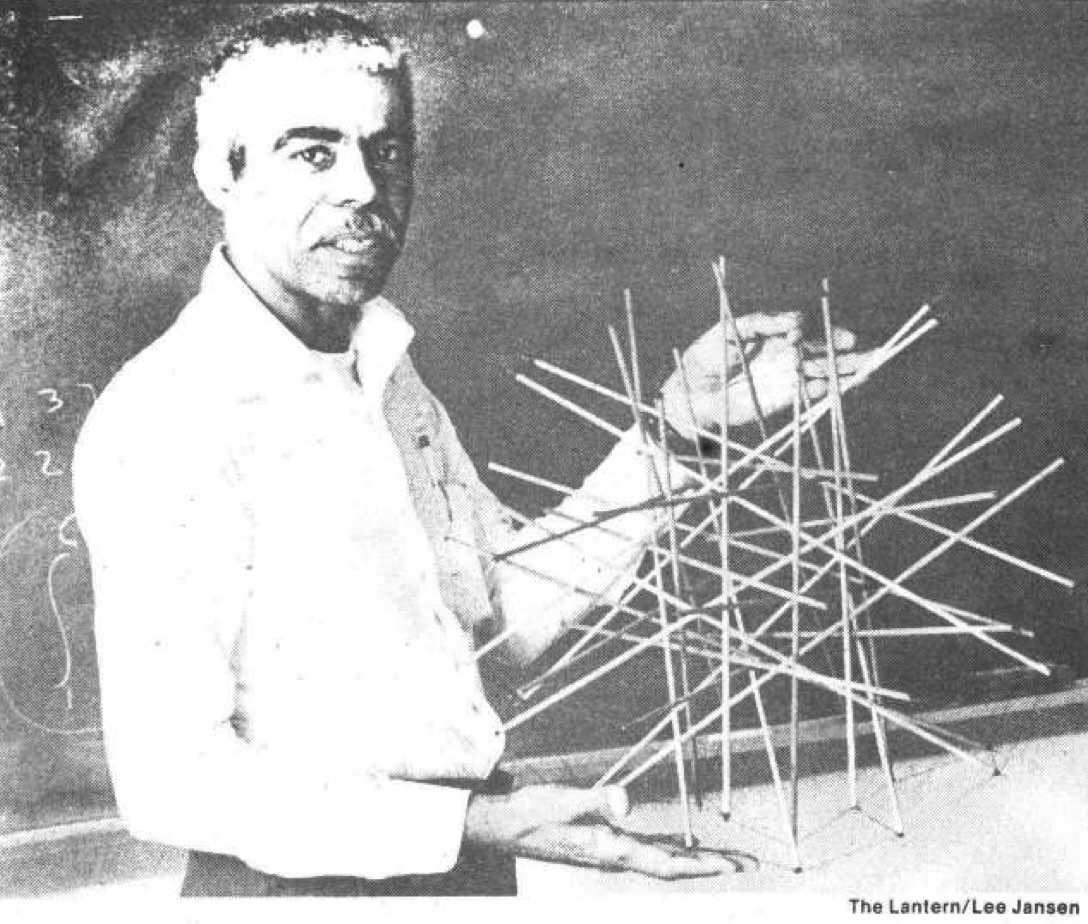}
    \end{center}
    \caption{William McWorter, Jr. in 1982 as an Associate Professor of Mathematics at The Ohio State University. (Photo credit from the Lantern \cite{lantern1982nosticks}.)}
    \label{fig:mcworter_portrait}
\end{figure}

This archive piqued my interest in McWorter. I wanted to find more, and we did. A black and white photo in the university newspaper in 1982 shows a picture of him smiling with a structure he had built, created with 30 sticks that never touch (Figure~\ref{fig:mcworter_portrait}). 
A textbook for an honors linear algebra course, describing a unique pedagogical approach to the common task of computing the characteristic polynomial of a square matrix with real coefficients. Research records showing a breadth of work in fractal theory. 

We had focused on context in terms of the social and political circumstances of a time period, but I also wanted to understand the mathematical context. 
Along my journey as a research mathematician, I've found without this context my understanding of certain concepts, fields, and even individual results feels fragmented, incomplete. And yet, many of the mathematicians I look up to seem to have an innate ability to connect theories, fields, and histories of results into a larger mathematical picture. Maybe stories can help us make those connections. This project was an opportunity for me to use storytelling not just to learn about mathematicians, but to learn mathematics. McWorter, in particular, gave me a clear place to begin. What is a fractal? Why did he care about them? What was happening in the mathematics of fractal geometry while he was working in the field? And what might his work tell us about the mathematical community he was part of?

\bigskip

\noindent\hdashrule{\linewidth}{0.5pt}{2pt 2pt}

\bigskip

\begin{figure}[htpb]
\begin{center}

\begin{tikzpicture}
    \draw
        (210:1.4)
        -- (90:1.4)
        -- (330:1.4)
        -- cycle;
    \node at (0,-1.8) {$n=0$};
\end{tikzpicture}
\hspace{0.5cm}
%
\begin{tikzpicture}
    \draw[decoration=Koch snowflake]
        decorate{
            (210:1.4)
            -- (90:1.4)
            -- (330:1.4)
            -- cycle
        };
    \node at (0,-1.8) {$n=1$};
\end{tikzpicture}
\hspace{0.5cm}
%
\begin{tikzpicture}
    \draw[decoration=Koch snowflake]
        decorate{
            decorate{
                (210:1.4)
                -- (90:1.4)
                -- (330:1.4)
                -- cycle
            }
        };
    \node at (0,-1.8) {$n=2$};
\end{tikzpicture}
\hspace{0.5cm}
%
\begin{tikzpicture}
    \draw[decoration=Koch snowflake]
        decorate{
            decorate{
                decorate{
                    (210:1.4)
                    -- (90:1.4)
                    -- (330:1.4)
                    -- cycle
                }
            }
        };
    \node at (0,-1.8) {$n=3$};
\end{tikzpicture}
\end{center}
\caption{Four iterations of the Koch snowflake}
\label{fig:koch_fractal}
\end{figure}
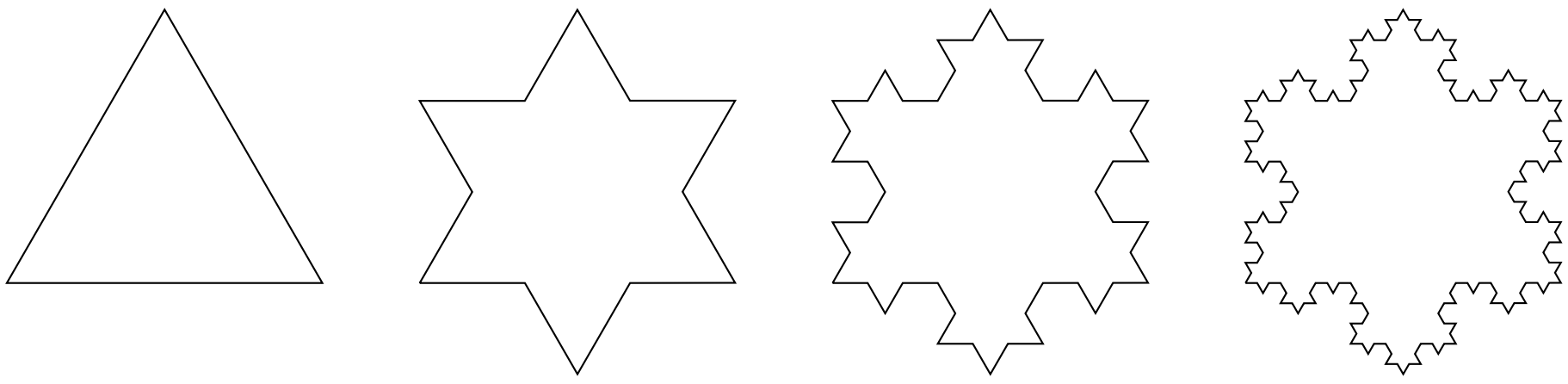

Fractals are geometric objects whose structure can persist across different scales. Some fractals exhibit self similarity, when you zoom you see patterns that resemble the larger object. Fractal like objects were being studied long before the word ``fractal'' was introduced. In the nineteenth century, Karl Weierstrass produced an example of a continuous function that is nowhere differentiable, challenging the familiar picture of a curve as something smooth. In 1904, Helge von Koch gave a geometric construction of what is now
called the Koch curve. It begins with a line segment, and each segment is iteratively replaced by four smaller segments, each of which is 1/3 the length of the original segment. Applying this construction to all three sides of an equilateral triangle produces the Koch snowflake (Figure~\ref{fig:koch_fractal}). Repeating this process produces increasingly intricate boundary, while retaining the same basic pattern at every scale. 

The Koch curve was the first time I saw that dimension does not always have to be an integer. I am accustomed to thinking of a line as one dimensional and a region of the plane as two dimensional. The Koch curve consists of four copies of itself, each scaled by a factor of $1/3$. For a self similar object consisting of $N$ copies scaled by a factor $r$, the \emph{similarity dimension} $d$ is determined by,
$$
Nr^d=1.
$$
For the Koch curve, this gives,
$$
4\left(\frac{1}{3}\right)^d=1,
$$
and therefore the dimension $d$ is,
$$
d=\frac{\log 4}{\log 3}\approx1.2619.
$$
The resulting dimension lies between $1$ and $2$. The curve has a richer
self similar structure than an ordinary line, but it does not fill a
two dimensional region.

These ideas developed across several areas of mathematics before there was a field called fractal geometry. Work in analysis, set theory, topology, and complex dynamics produced many of the objects and ideas that would later be brought together under the language of fractals by Benoît Mandelbrot. While working at IBM's Thomas J. Watson Research Center, he used computers to investigate and visualize irregular and self similar structures, connecting mathematical ideas to phenomena such as coastlines and other naturally occurring forms of roughness. He introduced the term \emph{fractal} in 1975 and helped establish fractal geometry as a mathematical research area.

\begin{figure}[htpb]
    \begin{center}
    \includegraphics[width=\textwidth]{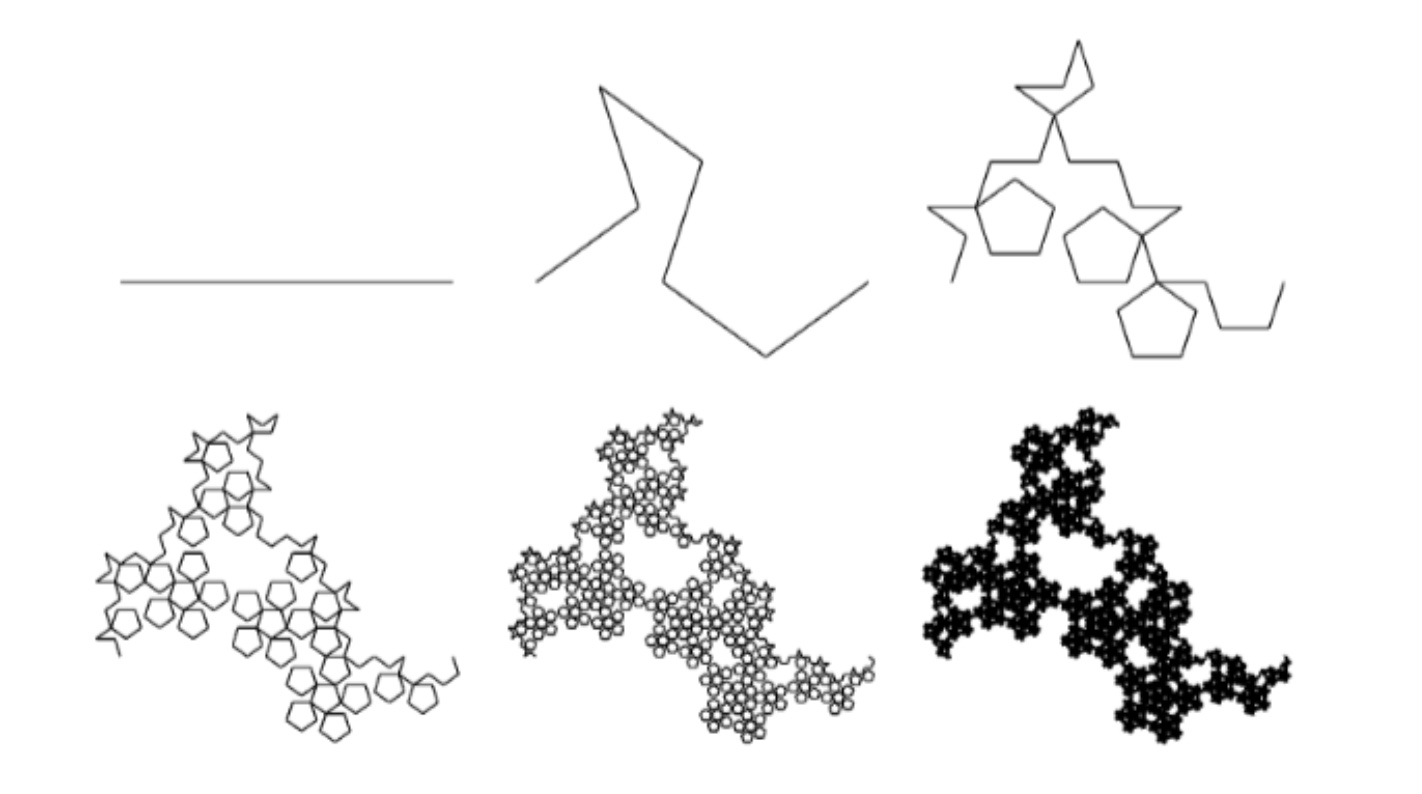}
    \end{center}
    \caption{Six iterations of the McWorter Pentigree fractal from Edgar's \emph{Measure, Topology, and Fractal Geometry} \cite{edgar1990measure} .}
    \label{fig:pentigree_fractal}
\end{figure}

This history changed the way I understood McWorter's work. By the 1980s,
mathematicians were studying self-similarity, dimension, iteration, and
irregular geometric structure with the aid of computers. McWorter's
tenure at Ohio State overlapped with this period of rapid development.
His work included experiments with recursive constructions generated by
computer programs. In 1987, he published an article entitled
``Creating Fractals'' \cite{mcworter1987creating}, which introduced and
detailed the Pentigree, a self similar fractal curve with five fold symmetry that has since been named after him (Figure~\ref{fig:pentigree_fractal}).

The Pentigree begins with a line segment, and from this a six segment motif is constructed using rotations. The resulting scale factor is:

$$
r=\frac{1}{1+2\cos(36^\circ)}
=\frac{3-\sqrt{5}}{2}
\approx0.381966.
$$

\noindent Each segment is then replaced recursively by a scaled and rotated copy
of the six segment motif. The same idea we used for the Koch curve applies to the Pentigree. Each of the six segments generates a smaller copy of the entire curve, scaled by $r$ so its similarity dimension $d$ satisfies
$
6r^d=1.
$
Thus,

$$
d=\frac{\log(1/6)}{\log(r)}
\approx1.86172.
$$

One way mathematicians formalize these constructions is through an
\emph{iterated function system}, or IFS, a collection of transformations that are repeatedly applied to an object. In the Pentigree, the six transformations each scale, rotate, and translate a copy of the entire curve. A small collection of geometric instructions can therefore generate increasingly complicated structure. McWorter's 1987 article described how recursive algorithms could be used to generate such fractal constructions with a computer. Around the same time, his Ohio State colleague Gerald Edgar wrote \emph{Measure, Topology, and Fractal Geometry} \cite{edgar1990measure}, a text that grew from notes associated with a summer course he taught for the Ross Mathematics Program in 1988. Chapter 1 introduces McWorter's Pentigree. 

A box containing family items, among them a copy of a mathematical article, led me into the world and history of fractals in the United States and their intersections with computing. It also gave me an appreciation for what productivity could look like when working in a growing research area, with like minded colleagues and graduate students down the hall.\footnote{Several students in the 1990s had theses dealing explicitly with fractal theory. This suggests Edgar's group contributed to an intellectual fractal theory community within the department during this period \cite{mathgenealogy11345}.}. McWorter's story was beginning to take shape across all three dimensions of our framework: his family and career, his mathematics, and the mathematical and historical context around him. This introduced another challenge for our project. How could we help someone else encounter that mathematics?

\bigskip

\noindent\hdashrule{\linewidth}{0.5pt}{2pt 2pt}

\bigskip

We wanted to introduce others to the stories of Black mathematicians from Ohio State. How do you introduce someone to a mathematician they have never heard of, whose mathematics they've never encountered? How do you help students see research mathematics as something people actually do? We decided to introduce their math through play. We designed workshops for middle school students featuring three parts: focusing on humanity, mathematics, and context.

\begin{figure}[htpb]
    \begin{center}
    \includegraphics[width=\textwidth]{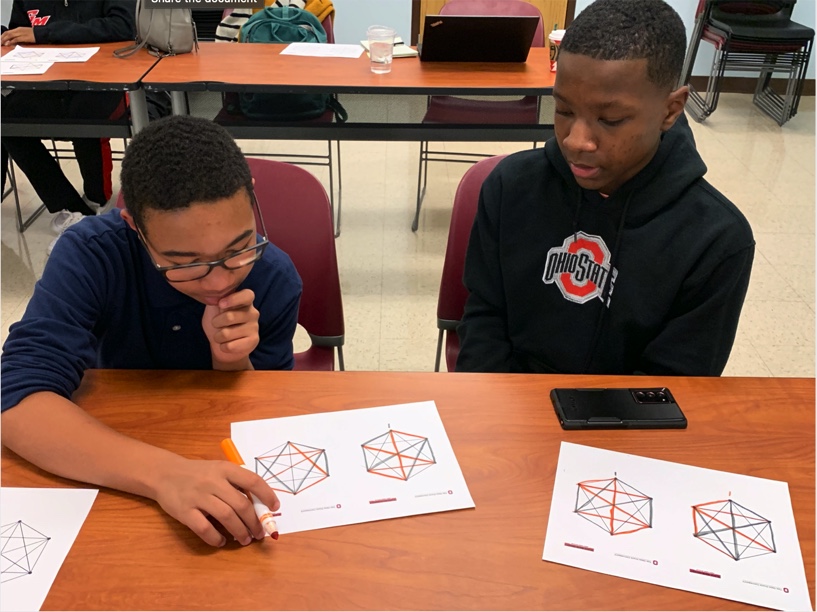}
    \end{center}
    \caption{Middle School Students work on Matroid Puzzle Game inspired by the Mathematics of Carolyn Mahoney in a Hidden Figures Revealed Workshop.}
    \label{fig:hid_fig_wkshp}
\end{figure}

\begin{figure}[htpb]
    \begin{center}
    \includegraphics[width=\textwidth]{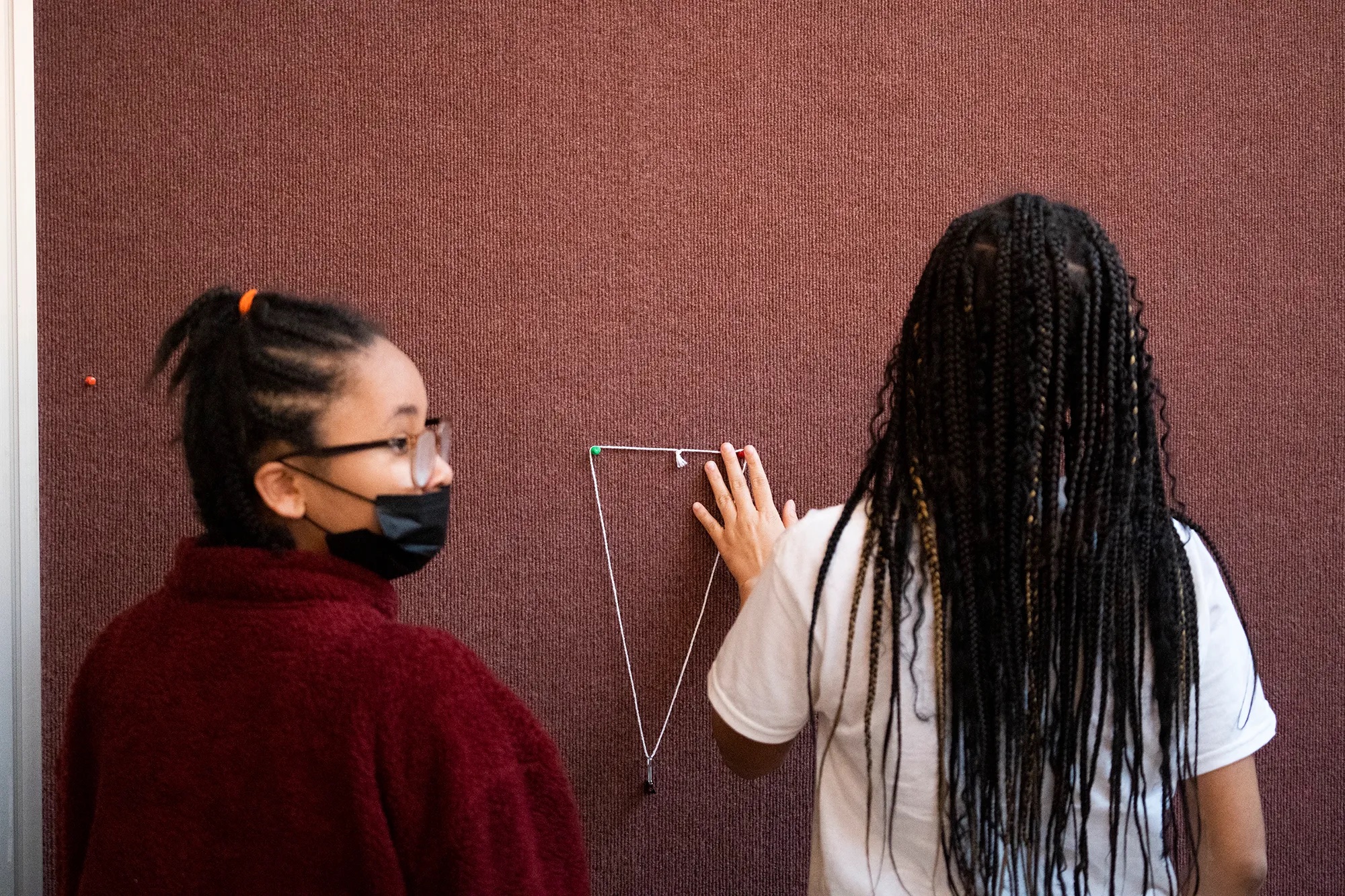}
    \end{center}
    \caption{Middle School Students work on a puzzle hanging activity inspired by the mathematics of Guy Hogan in a Hidden Figures Revealed Workshop. (Photo credit from  \cite{hendrix2023meet}.)}
    \label{fig:hid_fig_wkshp_hogan}
\end{figure}

Carolyn Mahoney worked in matroid theory. I learned with my students that a matroid is a mathematical structure that abstracts the idea of independence from linear algebra. Formally, a matroid $M$ is a pair $M=(E,\mathcal{I})$, where $E$ is the ground set of $M$ and $\mathcal{I}$ is a collection of subsets of $E$ called independent sets. Given an undirected graph $G=(V,E)$, we can associate to it a graphic matroid $M(G)=(E,\mathcal{I})$, where $\mathcal{I}$ is the set of subsets of $E$ that are acyclic. We made a modified version of “Dots and Boxes” inspired by this idea (Figure~\ref{fig:hid_fig_wkshp}).\footnote{In the game, students took turns drawing lines on a grid to enclose boxes and claim the enclosed area. They lost a point if they created a cycle. Mathematically, creating a cycle meant the corresponding set of edges was no longer independent. Of course, we did not lead with this language while they were playing!} Guy Hogan studied group theory, so we designed a puzzle hanging activity that introduced the idea of operations through a physical sequence of moves (Figure~\ref{fig:hid_fig_wkshp_hogan}). A plus meant moving a weight on a string to the right over two thumb tacks, a minus meant moving the string to the left. What did plus plus minus plus do? Was there a way to create a sequence that would create the weighted string to fall with the removal of just one pin? It seemed impossible.

We started with humanity and spent about 30 minutes with the students before introducing the mathematics. We asked them what they thought mathematicians did, what they knew about research. Had they ever heard of the Republic of Panama? That's where Guy Hogan grew up. He went to an HBCU, Talladega College, in Alabama, had they ever heard of HBCUs? Did they know anyone who left the state to go to college? For
Mahoney, we talked about her work developing mathematics curricula in schools and asked what, if anything, students would change about mathematics education in their own schools. Hogan eventually left mathematics for a legal career, so we asked whether their favorite subjects had ever changed. 

Why do people care about operations? In the puzzle game we composed moves that behaved like operations, but there were no numbers involved. Mathematicians like Guy Hogan study operations abstractly in subjects like group theory. We talked about the movement of robots and how rotations can be viewed as operations that can be composed. We talked about video games and how they use mathematics to help the characters move across the screen. For the matroid activity, we returned to graphs. We had students write their favorite musical artists as nodes on a board, then they decided to put edges between them if they were the same type of artist. We discussed how this could be a helpful way to organize connections. We talked about social media.

Our archival work gradually developed into a framework for storytelling that shaped the narratives we created about Black math alumni. The principles of humanity, mathematics, and context guided our research and ultimately became a way to facilitate mathematical discovery, both for middle school students and for ourselves. The workshops made the framework operational. The goal was not simply to tell students that Black mathematicians had done interesting work. We want to show them new mathematical ideas, learn about the people whose ideas inspired the activity, and understand something about the worlds in which they lived. This gave a proof of concept that our framework could do more than shape a story, it could also shape how someone encounters mathematics. 

 I hope you reflect on these principles of mathematical storytelling, \emph{mathematics, humanity, and context}, and think about how they could work in practice. From now on, listen closely to the mathematical stories around you and how they are told. How often do they detail history? How often do they make the mathematics itself visible? How often do they allow you to encounter the person who did the mathematics? If one of these pieces is missing, ask what becomes harder to see. We can, and should, do more. We can open up our mathematical community to ourselves and to the public.


\letterbreak


\begin{quote}
\emph{``I am what time, circumstance, history, have made of me, certainly, but I am also, much more than that. So are we all."}\\
    
\hfill -- James A. Baldwin \cite{baldwin1984notes} --
\end{quote}

If you haven't heard it already, there's a pervasive toolbox metaphor in mathematics. The idea is that we all have a box that we fill over time with different techniques and methods. Each class you take, each research experience, each conversation with another mathematician is an opportunity to add something new. I hope you find ways to travel, to attend conferences and lectures and seminars, to build your understanding of theory, to learn about applications that interest you, and to spend time thinking deeply about problems that linger in your head for weeks and months. Each of these experiences connects you to more mathematicians, to more ideas, and gives you more to add to your box. Even when it feels like you don't know enough, stop and peek. There are so many pieces waiting to be connected. Have conversations. Bridge the gaps. You have the capacity for this. Let yourself grow.

The archival research I did through Hidden Figures Revealed has inspired my present work to do a Census of Black mathematicians who earned PhDs in the mathematical sciences across the United States. The NSF Ascending Postdoc helped me add tools related to data science to my mathematical toolbox. Now I view information from these historical studies as data. I can consider degrees earned and work experiences, patterns of matriculation, institutional pathways, and research networks from quantitative instead of qualitative perspectives and use applied mathematical tools to study the mathematical community. I could not have imagined this when I began my postdoc years ago. It did not occur to me that I could write my own mathematical story and build a research program that centers equity related questions concerning shape, structure, and society, part of the emerging field of quantitative justice.

Your mathematical life will probably not look the way you imagine it now. You will find questions you did not know you wanted to ask. You will learn mathematics you did not know existed. You will meet people who change the way you think. You will discover parts of yourself that you had not yet figured out how to bring with you into mathematics. I hope that you learn stories that connect you to mathematics, and pull you in and keep you here. I hope that if you want to keep going, you find pathways into the mathematical professions of your dreams. To do this you will have to keep adding tools as you go. The mathematical world you are entering is changing fast.\footnote{I am writing this as things like anti intellectualism, the demographic cliff in higher education, and the rapid incorporation of AI into research pipelines create seismic shifts in what it means to do mathematics.} You will have to grow in ways that are unique to your generation. You will encounter new problems and new opportunities that I cannot anticipate.

You may struggle to reconcile your love for mathematics amid cultural tension, societal unrest, and an endless list of things to learn. You may face the pang of disappointment that follows when you realize the toolboxes of many, when it comes to centering inclusion and diversity, are empty. You may have peers or colleagues who are increasingly frustrated at hearing sentiments like this expressed, who are tired of conversations about inclusion, who believe enough talent will supersede it all, that the math that draws you together exists independent from the politics of our society. Who don't want stories with context, and don't want stories with humanity. They believe the stories worth telling center mathematics and mathematics alone. They also do not want this romantic notion challenged, preferring the fairytale, a merit based utopia where the cream of the crop rise no matter what. You will have to decide how to navigate this. I hope you do not mistake their discomfort for evidence that the mathematical community is not for you. I hope you remember that your job is not to convince someone that you belong.

Storytelling is a tool. I want it to be featured in more mathematical toolboxes. I want you to consider adding it to yours. I don't mean that every mathematician needs to become a historian. I don't mean that every mathematics course needs to become a history course. And I certainly don't mean that storytelling is a solution to the complicated structural reasons that people enter, leave, or remain in mathematics. But I think there is a practice here that can be learned, adapted, and used when it makes sense. I think we can learn to tell stories that introduce people, not just math majors or graduate students or members of the professoriate or industry, but the general public, to mathematics rather than simply decorating mathematics with biographies. I think we can recover histories that help people see themselves as part of a mathematical community. I think we can ask who someone was, what mathematics they did, and what was happening around them, and let all three answers matter.

Short, clean stories can dilute the complexity of history. Brown v. Board of Education outlawed `separate but equal' in public education in 1954, but Columbus, Ohio, continued to confront school segregation for decades, including through court ordered desegregation and busing beginning in 1979. Fractal geometry was not simply invented by Mandelbrot in 1975. There was a longer history of work on irregular and self similar objects before this theory was formalized. There were other classes of fractals being generated by computers during the time period that the Mandelbrot set was gaining popularity, including those attributed to William McWorter, Jr. The history of Black mathematicians at Ohio State cannot be reduced to a story of isolated ``firsts.'' In 2021, I was the only Black woman mathematician in the Department of Mathematics at The Ohio State University, yet between 1963 and 2020 there were 7 PhDs, 20 master's degrees, and 141 bachelor's degrees awarded to Black mathematicians. 

You may enter many rooms and realize you are the only person who looks like you. I  experience this, so did the Black mathematicians from Ohio State. If you find yourself standing at the edge of these rooms wondering whether there is a place for you inside, literally or figuratively, I hope you remember that there are more stories than the ones you have been told. That is perhaps the thing I most want you to take from this letter. Stories can make visible what has been hidden. Stories can transform a name into a person, a theorem into a history, and a career into something shaped not only by individual choices but by relationships, institutions, opportunities, and circumstances. As wonderful as they can be, stories can also flatten. Please resist this. Use them as a tool to push against prevailing, but reductive, narratives, especially the ones that hand wave why entire demographics are less represented in our community. 

That is what I did in some sense. I challenged a narrative that there had been no one like me in the math department at Ohio State before. I found Rada Higgins, then I found hundreds of other names. The names led to archives, the archives led to mathematics, and the mathematics led to questions about mathematicians and the worlds in which they worked. I started looking for evidence that I was not alone, and in seeking to tell their stories, I developed a framework centering humanity, context, and mathematics as a way to tell stories about mathematicians in general.

Your mathematical story is still being written. Lean in as you encounter stories around you. Think deeply about your own. Ask the right questions. Who are you? Where are you from? What mathematics do you do? What narratives have you told yourself about what is possible for your future? What assumptions have you made about your place in the mathematical community, and what have you been told? What becomes possible to imagine when you consider the communities and context in which you learn and do your work? What needs to stay the same? What needs to change?

I hope Hidden Figures Revealed gives you something concrete to hold onto, a model of earnest efforts to more fully understand the history that shaped the mathematical community you are entering. So you might learn from those who came before you, understand the stories of Black mathematicians, and begin to have a more nuanced understanding of why there is a dearth of Black mathematicians. This project was also for anyone who loves mathematics, especially those whose stories are rarely centered. It is for anyone who has ever felt alone in mathematics, who has pushed past the fear that the weight of an entire community hinges on their success. This project was as much for me as it was for all of you. It was my way to reconcile my love of math with the rest of it. I hope you find your own way.

\bigskip

\noindent \normalfont{Mathematically yours},

\medskip

\noindent \Large{{\calligra A fellow mathematician}}




\bibliography{biblio}

\end{document}